\documentstyle[11pt,newlfont]{amsart}

\newcommand{\nc}{\newcommand}

\nc{\one}{\mbox{\bf 1}}
\nc{\invtensor}{\underset{\leftarrow}{\otimes}}
\nc{\rlarrows}{\begin{picture}(1,0.4)
                \put(0,-0.1){\makebox(1,0.2){$\leftarrow$}}
                \put(0,0.1){\makebox(1,0.2){$\ra$}}
              \end{picture}}
\nc{\rra}{\begin{picture}(1,0.4)
                \put(0,-0.1){\makebox(1,0.2){$\lra$}}
                \put(0,0.1){\makebox(1,0.2){$\lra$}}
              \end{picture}}
\nc{\Left}{\mathbf L}  
\nc{\Right}{\mathbf R} 

\nc{\gr}{\operatorname{gr}}
\nc{\Ho}{\operatorname{Ho}}
\nc{\alt}{\operatorname{alt}}
\nc{\Sym}{\operatorname{Sym}}
\nc{\sym}{\operatorname{sym}}
\nc{\id}{\operatorname{id}}
\nc{\Der}{\operatorname{Der}}
\nc{\im}{\operatorname{Im}}
\nc{\Ker}{\operatorname{Ker}}
\nc{\Col}{\operatorname{Col}}
\nc{\ter}{\operatorname{ter}}
\nc{\intl}{\operatorname{int}}
\nc{\out}{\operatorname{out}}
\nc{\irr}{\operatorname{irr}}
\nc{\Tor}{\operatorname{Tor}}
\nc{\WE}{\frak W}
\nc{\FIB}{\frak F}
\nc{\COF}{\frak C}
\nc{\Sh}{\operatorname{Sh}}
\nc{\Hom}{\operatorname{Hom}}
\nc{\uhom}{\operatorname{\ul{ hom}}}
\nc{\End}{\operatorname{End}}
\nc{\holim}{\operatorname{holim}}
\nc{\dirlim}{\underset{\rightarrow}{\lim}\,}
\nc{\invlim}{\underset{\leftarrow}{\lim}\,}
\nc{\CB}{\operatorname{\bf CB}}
\nc{\Op}{\operatorname{Op}}
\nc{\com}{\operatorname{co}}
\nc{\Tot}{\operatorname{Tot}}
\nc{\Th}{\operatorname{Th}}
\nc{\Cech}{\check{C}}
\nc{\Spec}{\operatorname{Spec}}
\nc{\Spf}{\operatorname{Spf}}
\nc{\MC}{\operatorname{MC}}
\nc{\U}{\operatorname{U}}
\nc{\Diff}{{\cal D}\mbox{\em iff}}
\nc{\Mor} {{\cal M}or}
\nc{\Ob}{\operatorname{Ob}}
\nc{\cone}{\widehat}
\nc{\Coder}{\operatorname{Coder}}
\nc{\Lie}{\operatorname{Lie}}
\nc{\dbar}{\bar{\partial}}
\nc{\pr}{\operatorname{pr}}
\nc{\diag}{\operatorname{diag}}

\nc{\Mod}{{\mathtt{mod}}}       
\nc{\Ab}{{\mathtt {Ab}}}          
\nc{\dgl}{{\mathtt{dglie}}}
\nc{\dga}{{\mathtt{dgalg}}}
\nc{\art}{{\mathtt {art}}}
\nc{\simpl}{\Delta^0\Ens}
\nc{\Kan}{{\mathtt {Kan}}}

\nc{\Grp}{{\mathtt {Grp}}}
\nc{\dgc}{{\mathtt {cdga}}}
\nc{\Ens}{{\mathtt {Ens}}}

\nc{\pa}{\partial}

\nc{\CA}{\cal A}
\nc{\CF}{\cal F}
\nc{\CG}{\cal D}
\nc{\CI}{\cal I}
\nc{\CJ}{\cal J}
\nc{\CO}{\cal O}
\nc{\CU}{\cal U}
\nc{\CC}{\cal C}
\nc{\CDD}{\cal D}
\nc{\CL}{\cal L}
\nc{\CM}{\cal M}
\nc{\CS}{\cal S}
\nc{\CT}{\cal T}

\nc{\fg}{\frak g}
\nc{\fk}{\frak k}
\nc{\fh}{\frak h}
\nc{\fm}{\frak m}
\nc{\fn}{\frak n}
\nc{\fS}{\frak S}

\nc{\nen}{\newenvironment}
\nc{\ol}{\overline}
\nc{\ul}{\underline}
\nc{\ra}{\rightarrow}
\nc{\lra}{\longrightarrow}
\nc{\lla}{\longleftarrow}
\nc{\Lra}{\Longrightarrow}
\nc{\Lla}{\Longleftarrow}
\nc{\Llra}{\Longleftrightarrow}
\nc{\hra}{\hookrightarrow}
\nc{\iso}{\overset{~}{\lra}}

\nc{\Thm}[1]{Theorem~\ref{#1}}
\nc{\Prop}[1]{Proposition~\ref{#1}}
\nc{\Lem}[1]{Lemma~\ref{#1}}
\nc{\Cor}[1]{Corollary~\ref{#1}}
\nc{\Conj}[1]{Conjecture~\ref{#1}}
\nc{\Claim}[1]{Claim~\ref{#1}}
\nc{\Defn}[1]{Definition~\ref{#1}}
\nc{\Exa}[1]{Example~\ref{#1}}
\nc{\Rem}[1]{Remark~\ref{#1}}
\nc{\Note}[1]{Note~\ref{#1}}

\nen{thm}[1]{\label{#1}{\bf Theorem.\ } \em}{}
\nen{prop}[1]{\label{#1}{\bf Proposition.\ } \em}{}
\nen{lem}[1]{\label{#1}{\bf Lemma.\ } \em}{}
\nen{cor}[1]{\label{#1}{\bf Corollary.\ } \em}{}
\nen{conj}[1]{\label{#1}{\bf Conjecture.\ } \em}{}
\nen{claim}[1]{\label{#1}{\bf Claim.\ } \em}{}

\nen{defn}[1]{\label{#1}{\bf Definition.\ } }{}
\nen{exa}[1]{\label{#1}{\bf Example.\ } }{}

\nen{rem}[1]{\label{#1}{\em Remark.\ } }{}
\nen{note}[1]{\label{#1}{\em Note.\ } }{}
\nen{exer}[1]{\label{#1}{\em Exercise.\ } }{}

\begin{document}

\title[]{Descent of Deligne groupoids}
\author{Vladimir Hinich}
\address{Dept. of Mathematics and Computer Science, University of Haifa,
Mount Carmel, Haifa 31905 Israel}


\maketitle

\begin{abstract}
To any non-negatively graded dg Lie algebra $\fg$ over a field $k$ of
characteristic zero  we assign a functor $\Sigma_{\fg}:\art/k\to\Kan$ from the
category of commutative local artinian $k$-algebras with the residue field $k$
to the category of Kan simplicial sets. There is a natural homotopy equivalence
between $\Sigma_{\fg}$ and the Deligne groupoid corresponding to $\fg$.

The main result of the paper claims that the functor $\Sigma$ commutes 
up to homotopy with the "total space" functors which assign a dg Lie 
algebra to a cosimplicial dg Lie algebra and a simplicial set to a 
cosimplicial simplicial set. This proves a conjecture of Schechtman
~\cite{s,s1,hs} which implies that if a deformation problem 
is described ``locally'' by a sheaf of dg Lie algebras $\fg$ on a 
topological space $X$ then the global deformation problem is described 
by the homotopy Lie algebra $\Right\Gamma(X,\fg)$. 
\end{abstract}

\section{Introduction}

\subsection{}
\label{i1}
Let $\fg$ be a dg Lie algebra over a field $k$ of characteristic zero
concentrated in non-negative degrees. The algebra $\fg$ defines a functor
$$ \CC_{\fg}:\art/k\to\Grp$$
from the category of local artinian $k$-algebras with the residue field $k$
to the category of groupoids --- see~\ref{del-gr} or~\cite{gm1}, sect.~2.

A common belief is that any "reasonable" formal deformation problem
can be described by the functor $\CC_{\fg}$ where $\fg$ is an
appropriate "Lie algebra of infinitesimal automorphisms". This would 
imply, for instance, that if $H^0(\fg)=0$ (i.e. if the automorphism group
of the deformed object is discrete) then the completion of the local ring
of a moduli space at a given point is isomorphic to the completion of
the $0$-th cohomology group of $\fg$.

If we are dealing with the deformations of  algebraic structures
(associative, commutative or Lie algebras or so), the
Lie algebra $\fg$ is just the standard complex calculating the cohomology
of the appropriate type.

In this paper we prove the following claim conjectured by V.~Schechtman in
~\cite{s,s1,hs}. It allows one to construct
a dg Lie algebra which governs various formal deformations in the non-affine
case.

{\em Let $X$ be a topological space and let $\fg$ be a sheaf of
dg Lie $k$-algebras. Let $\CU=\{U_i\}$ be an locally finite open covering
of $X$, $\CG_i=\CC_{\Gamma(U_i,\fg)}$ and let $\CG$ be the groupoid  of
"descent data" for the collection $\CG_i$ (see~\ref{tot:simpl},{\em2}). 
Then $\CG$ is naturally
equivalent to the groupoid $\CC_L$ where $L$ is a dg Lie algebra
representing the Cech complex $\Cech(\CU,\fg)$.}

This result implies, in particular, \Cor{generalization}
which claims that Theorem 8.3 of~\cite{hdtc}
remains valid without the assumption of formal smoothness.

\subsection{Structure of the Sections}
In Section~\ref{sigma} we define the {\em contents} $\Sigma(\fg)$ of
a nilpotent dg Lie algebra $\fg$. It is a Kan simplicial set homotopically
equivalent to the Deligne groupoid $\CC(\fg)$. 
In Section~\ref{tot} we recall the definition of the total space functor 
in different categories. We also present~\Prop{criterion} giving a sufficient
condition for a map of bisimplicial dg algebras to be an acyclic
fibration in the sense of~\ref{bisim-alg}. 
Now the claim~\ref{i1} can be interpreted as the 
commutativity of $\Sigma$ with the functors $\Tot$. The main 
\Thm{main} is proven in Section~\ref{thm}. In the last Section~\ref{app}
 we deduce from~\Thm{main} an application to formal deformation theory.

The idea that the generalization of the result of~\cite{hdtc} to the
non-smooth case should follow from a descent property for Deligne
groupoids belongs to V.~Schechtman. I am very grateful to him for
helpful discussions on the subject. \Prop{kan} claiming that the content
$\Sigma(\fg)$ and the Deligne groupoid $\CC(\fg)$ are homotopy equivalent
is consonant to the Main Homotopy Theorem of Schlessinger--Stasheff, cf.~
\cite{ss}. 
I am greatly indebted to J.~Stasheff who read the first draft of the manuscript
and made some important remarks. 

\subsection{Notations}

Throughout this paper $k$ is a fixed field of characteristic zero.

$\art/k$ denotes the category of commutative local artinian $k$-algebras
having the residue field $k$.

$\dgl(k)$ (resp., $\dgc(k)$) is the category of non-negatively graded dg Lie 
(resp., commutative) algebras over $k$.

$\Delta$ is the category of ordered sets $[n]=\{0,\ldots,n\},\ n\geq 0$
and monotone maps; $\simpl$ is the category of simplicial sets;
$\Delta^n\in\simpl$ are the standard $n$-simplices.

$\Kan\subseteq\simpl$ is the full subcategory of Kan simplicial sets.

$\Ab$ is the category of abelian groups. 
$C(\CA)$ (sometimes $C(R)$) is the category of complexes over an abelian
category $\CA$ (over the category of $R$-modules). $C^{\geq 0}$ denotes the
full subcategory of non-negatively graded complexes.

\section{Contents of a nilpotent Lie dg algebra}
\label{sigma}

\subsection{}

For any $n\geq 0$ denote by $\Omega_n$ the $k$-algebra of polynomial
differential forms on the standard $n$-simplex $\Delta^n$ --- 
see~\cite{bg}. 

One has 
$$\Omega_n=k[t_0,\ldots,t_n,dt_0,\ldots,dt_n]/(\sum t_i-1,\sum dt_i).$$

The algebras $\Omega_n$ form a simplicial commutative dg algebra:
a map $u:[p]\to[q]$ induces the map $\Omega(u):\Omega_q\to\Omega_p$
defined by the formula $\Omega(u)(t_i)=\sum_{u(j)=i}t_j$.

If $\fg$ is a dg Lie $k$-algebra and $A$ is a commutative dg $k$-algebra
then the tensor product $A\otimes\fg$ is also a dg Lie $k$-algebra. Thus, 
any dg Lie algebra $\fg$ gives rise to a simplicial dg Lie algebra
$$\fg_{\bullet}=\{\fg_n=\Omega_n\otimes\fg\}_{n\geq 0}.$$

For any dg Lie algebra $\fg$ denote by $\MC(\fg)$ the set of elements
$x\in\fg^1$ satisfying the Maurer-Cartan equation:
$$ dx+\frac{1}{2}[x,x]=0.$$

\subsubsection{}
\begin{defn}{contents}
Let $\fg\in\dgl(k)$ be nilpotent. Its {\em contents} 
$\Sigma(\fg)\in\simpl$
is defined as 
$$\Sigma(\fg)=\MC(\fg_{\bullet}).$$
\end{defn}

\subsubsection{} Recall (see~\cite{bg}) that the collection of commutative dg algebras
$\Omega_n$ defines a contravariant functor
$$ \Omega:\simpl\to\dgc(k)$$
so that $\Omega(\Delta^n)=\Omega_n$ and $\Omega$ carries  direct limits 
in $\simpl$ to inverse limits.

\begin{lem}{repr}
Let $S\in\simpl$. There is a natural map
$$\MC(\Omega(S)\otimes\fg)\to\Hom(S,\Sigma(\fg))$$
which is bijective provided $S$ is finite (i.e., has a finite 
number of non-degenerate simplices).
\end{lem}
\begin{pf} This is because tensoring by $\fg$ commutes with finite
limits --- compare to~\cite{bg}, 5.2.
\end{pf}

\subsubsection{}
\begin{defn}{AF}
A map $f:\fg\to\fh$ of nilpotent algebras in $\dgl(k)$ will be called an
{\em acyclic fibration} if it is surjective and induces a quasi-isomorphism
of the corresponding lower central series.
\end{defn}

\subsubsection{}
\begin{lem}{AFsur}
Let $f:\fg\to\fh$ be an acyclic fibration of nilpotent dg Lie algebras in
$\dgl(k)$. Then the induced map $\MC(f):\MC(\fg)\to\MC(\fh)$ is 
surjective.
\end{lem}
\begin{pf}
Induction by the nilpotence degree of $\fg$ --- similarly 
to~\cite{gm1},~Th. 2.4.
\end{pf}

\subsubsection{}
\begin{lem}{A+L}Let $f:A\to B$ be a surjective map in $\dgc(k)$ and
$g:\fg\to\fh$ be a surjective map in $\dgl(k)$. Then the map
$$ A\otimes\fg\to (A\otimes\fh)\times_{(B\otimes\fh)}(B\otimes\fg)$$
is an acyclic fibration provided either
(a) $f$ is quasi-isomorphism or (b) $g$ is acyclic fibration.
\end{lem}
\begin{pf}
Since for any commutative dg algebra $A$ (with 1)
the functor $A\otimes\_$ transforms the lower central series of $\fg$
into the lower central series of $A\otimes\fg$, it suffices to check
that the above map is a surjective quasi-isomorphism.
 This is fairly standard.
\end{pf}

\subsubsection{}
\begin{prop}{AFaf}
Let $f$ be as above. Then $\Sigma(f):\Sigma(\fg)\to\Sigma(\fh)$ is an acyclic
fibration of simplicial sets.
\end{prop}
\begin{pf}
Lemmas~\ref{repr} and~\ref{AFsur} reduce the question to the following.
Let  $K\to L$ be an injective map of simplicial sets. Then one has to show
that the induced map of nilpotent Lie algebras
$$ \Omega(L)\otimes\fg\to\Omega(L)\otimes\fh\times_{\Omega(K)\otimes\fh}
\Omega(K)\otimes\fg$$
is an acyclic fibration. This follows from~\ref{A+L}(b).
\end{pf} 

\subsection{Deligne groupoid}
\label{del-gr}
Recall (cf.~\cite{gm1}) that for a nilpotent dg Lie algebra $\fg\in\dgl(k)$ 
defines the {\em Deligne groupoid} $\CC(\fg)$ as follows.

The Lie algebra $\fg^0$ acts on $\MC(\fg)$ by vector fields:
$$ \rho(y)(x)=dy+[x,y]\text{ for }y\in\fg^0, x\in\fg^1.$$
This defines the action of the nilpotent group $G=\exp(\fg^0)$ on the
set $\MC(\fg)$. Then the groupoid $\CC=\CC(\fg)$ is defined by the formulas
$$ \Ob\CC=\MC(\fg)$$
$$ \Hom_{\CC}(x,x')=\{g\in G|x'=g(x)\}.$$

\subsubsection{}
\begin{lem}{exp-1}
Let $\fg\in\dgl(k)$ be nilpotent. The natural map 
$\fg\to\fg_n=\Omega_n\otimes\fg$ induces an equivalence of groupoids 
$\CC(\fg)\to\CC(\fg_n)$.
\end{lem}
\begin{pf}
It suffices to check the claim when $n=1$. In this case an element
$z=x+dt\cdot y\in\fg_{1}=\fg[t,dt]$ with $x\in\fg^1[t],y\in\fg^0[t]$
satisfies MC iff $x(0)\in\MC(\fg)$ and $x$ satisfies the differential 
equation
$$ \dot{x}=dy+[x,y].$$
This  means precisely that $x=g(x(0))$ where $g\in G_1=\exp(\fg_1)$ is given
by the differential equation
$$ \dot{g}=g(y)$$
with the initial condition $g(0)=1$.
\end{pf}

\subsubsection{Explicit description of $\Sigma(\fg)$}
The notations are as in~\ref{del-gr}.

For any $n\geq 0$ let $G_n=\exp(\fg^0_n)$ be the group of polynomial maps
from the standard $n$-simplex $\Delta_n$ to the group $G$. The
collection $G_{\bullet}=\{G_n\}$ forms a simplicial group.  Right
multiplication defines on $G_{\bullet}$  a right $G$-action.

\begin{prop}{explicit}
There is a natural bijection 
$$G_{\bullet}\times^G\MC(\fg)\to\Sigma(\fg)$$
of simplicial sets. Here, as usual, $X\times^GY$ is the quotient 
of the cartesian product $X\times Y$ by the relation 
$$(xg,y)\sim (x,gy)\text{ for } x\in X,\ y\in Y,\ g\in G.$$
\end{prop}
\begin{pf} This immediately follows from~\Lem{exp-1}.
\end{pf}

Let $N\CC(\fg)\in\simpl$ be the nerve of $\CC(\fg)$.

Define the map $\tau:\Sigma(\fg)\to N\CC(\fg)$ by the formula
$$\tau(g,x)=(g_0(x),g_1g_0^{-1},\ldots,g_ng_{n-1}^{-1})$$
where $g\in G_n,x\in\MC(\fg), g_i=v_i(g)$ with $v_i$ being the 
$i$-th vertex.

\Prop{explicit} implies the following
\subsubsection{}
\begin{prop}{kan}
$\Sigma(\fg)$ is a Kan simplicial set. The map $\tau$ is an
acyclic fibration identifying $\CC(\fg)$ with the Poincar\'{e} groupoid
of $\Sigma(\fg)$.

More generally, if $f:\fg\to\fh$ is a surjective map of nilpotent
dg Lie algebras, then the induced map
$$ \Sigma(\fg)\to\CC(\fg)\times_{\CC(\fh)}\Sigma(\fh)$$
is an acyclic fibration.
\end{prop}
\begin{pf} The question reduces to the following. Given a pair of
polynomial maps $\alpha:\partial\Delta^n\to\fg^0,\ 
\beta:\Delta^n\to\fh^0$ satisfying $f\alpha=\beta|_{\partial\Delta^n}$
find a map $\gamma:\Delta^n\to\fg^0$ such that 
$\alpha=\gamma|_{\partial\Delta^n},\ \beta=f\gamma$. This is always
possible sincec the canonical map $\Omega_n\to\Omega(\partial\Delta^n)$
is surjective.
\end{pf}
\subsubsection{}
\begin{rem}{}
If one does not require $\fg$ to be non-negatively graded, its
contents is still a Kan  simplicial set. In this case it can
probably be considered as a generalization of the notion of
Deligne groupoid.
\end{rem}

\section{"Total space" functor and $\CM$-simplicial sets}
\label{tot}

\subsection{Generalities}
\subsubsection{The category $\CM$}
\label{catM}
Here and below $\CM$ denotes the following category of morphisms of
$\Delta$: 

The objects of $\CM$ are morphisms $[p]\to [q]$ in $\Delta$. 
A morphism from $[p]\to [q]$ to $[p']\to [q']$ is given by a 
commutative diagram
$$
\begin{array}{ccc}
[p] & {\lra} & [q] \\
{{\scriptstyle\alpha}\uparrow} &  & {{\scriptstyle\beta}\downarrow} \\
{[p']} & {\lra} & {[q']} \\
\end{array}
$$
The morphism in $\CM$ corresponding to $\alpha=\id,\beta=\sigma^i$,
is denoted by $\sigma^i$; the one corresponding to
$\alpha=\id,\beta=\partial^i$, is denoted by $\partial^i$.
"Dually", the morphism corresponding to $\alpha=\partial^i,\beta=\id$,
is denoted by $d_i$ and the one with $\alpha=\sigma^i,\beta=\id$, is
denoted by $s_i$.

\subsubsection{Total space}
Let $\CC$ be a simplicial category having inverse limits and functorial
{\em function objects} $\uhom(S,X)\in\CC$ for $X\in\CC,\ S\in\simpl$
--- see~\cite{bk}, IX.4.5 and the examples below. The total space $\Tot(X)$ 
of a cosimplicial object $X\in\Delta\CC$ is defined by the formula

$$ \Tot(X)=\invlim_{\phi\in\CM}\uhom(\Delta^p,X^q)$$
where $\phi:[p]\to [q]$.

\subsection{Examples} We will use three instances
of the described construction.

\subsubsection{}
\label{tot:simpl} Let $\CC=\simpl$. Then the above definition coincides 
with the standard one given in~{\em loc. cit.}, XI.3.

Let $G\in\Delta\Grp$ be a (strict) cosimplicial groupoid.
We will consider $\Grp$ as a full subcategory of $\simpl$,
so a simplicial set $\Tot(G)$ is defined.

\begin{lem}{tot-gr}
$T=\Tot(G)$ is  a groupoid. The objects of $T$ are collections 
$\{a\in\Ob G^0, \theta:\partial^1(a)\overset{\sim}{\ra}\partial^0(a)\}$
with  $\theta$  satisfying the cocycle condition:
$$ \sigma^0(\theta)=\id(a);
\ \partial^1\theta=\partial^0\theta\circ\partial^2\theta.$$

A morphism in $T$ from $\{a,\theta\}$ to 
$\{b,\theta'\}$  is a morphism $a\to b$ compatible with $\theta,\theta'$.

Thus, $\Tot(G)$ is ``the groupoid of descent data'' for $G$.
\end{lem}

\subsubsection{}
\label{tot:compl}
 Let $\CC=C(k)$ be the category of complexes over $k$.
For $S\in\simpl$ and $X\in C(k)$ the complex $\uhom(S,X)$ is defined
to be $\Hom(C_{\bullet}(S),X)$ where $C_{\bullet}$ is the complex
of normailized chains of $S$ with coefficients in $k$. The above
definition of the functor $\Tot$ coincides with the standard one.

\subsubsection{} 
\label{tot:lie}
Let $\CC=\dgl(k)$. For $S\in\simpl$ and $\fg\in\dgl(k)$
define $\uhom(S,\fg)=\Omega(S)\otimes\fg$. Then the functor
$\Tot:\Delta\dgl(k)\to \dgl(k)$ coincides with the Thom-Sullivan
functor described in~\cite{hdtc}, 5.2.4.
The De Rham theorem (see, e.g., {\em loc. cit.}, 5.2.8) shows that the 
functor $\Tot$ commutes up to homotopy with the forgetful functor
$ \#:\dgl(k)\to C(k).$

\subsection{$\CM$-simplicial sets}

In the sequel functors $X:\CM\to\simpl$ will be called 
$\CM$-simplicial sets. We now wish to find a sufficient condition for
a map $f:X\to Y$ of $\CM$-simplicial sets to induce an acyclic
fibration $\invlim f$ of the inverse limits. This will be an important
technical tool to prove the main theorem~\ref{main}.

\subsubsection{Matching spaces}
Fix $n\in\Bbb N$. Let $\partial^i:[n-1]\to[n]$, $i=0,\ldots,n$, be the 
standard face maps and let $\sigma^i:[n]\to[n-1]$, $i=0,\ldots,n-1$, be the
standard degeneracies.

The $n$-th matching space of a $\CM$-simplicial set $X$ is a 
simplicial subset $\mu_n(X)$ of the product
$$\prod_{i=0}^nX(\partial^i)\times\prod_{i=0}^{n-1}X(\sigma^i)$$ 
consisting of the collections 
$\left(x_i\in X(\partial^i),y^i\in X(\sigma^i)\right) $
satisfying the following three conditions:

($d$): $d_ix_j=d_{j-1}x_i\text{ for } i<j$.

($\sigma$): $\sigma^jy^i=\sigma^iy^{j+1}\text{ for } i\leq j$.

($d\sigma$): $\sigma^jx_i=d_iy^j\text{ for all } i,j$.

One has a canonical map $X(\id_n)\to\mu_n(X)$ which sends an
element $x\in X(\id_n)$ to the collection 
$(d_0x,\ldots,d_nx,\sigma^0x,\ldots,\sigma^{n-1}x)$.

\subsubsection{}
\begin{defn}{Maf}
A map $f:X\to Y$ of $\CM$-simplicial sets is called an {\em acyclic
fibration} if for any $n\in\Bbb N$ the commutative square
$$
\begin{array}{ccc}
X(\id_n) & {\lra} & Y(\id_n) \\
{\downarrow} &  & {\downarrow} \\
\mu_n(X) & {\lra} & \mu_n(Y) \\
\end{array}
$$
defines an acyclic fibration 
$$X(\id_n)\to Y(\id_n)\times_{\mu_n(Y)}\mu_n(X).$$
\end{defn}

\subsubsection{}
\begin{lem}{afM->ss}
Let $f:X\to Y$ be an acyclic fibration of $\CM$-simplicial sets.
Then the induced map of the corresponding inverse limits, $\invlim f$,
is an acyclic fibration.
\end{lem}
\begin{pf}  Let $\CM_{\leq n}$ be the full subcategory of $\CM$ 
consisting of morphisms $\alpha:[p]\to[q]$ with $p\leq n,q\leq n$.
Put $X(n)=\invlim X|_{\CM_{\leq n}}$. Then it is easy to see that
$X(n)=X(\id_n)\times_{\mu_n(X)}X(n-1)$. This immediately proves the
lemma.
\end{pf}

\subsection{Bisimplicial algebras}
\label{bisim-alg}
 Now we will present a source of various
$\CM$-simplicial sets in this paper.

\subsubsection{} Fix a cosimplicial nilpotent dg Lie algebra $\fg$.
Any bisimplicial commutative dg algebra  $A\in(\Delta^0)^2\dgc$
defines a $\CM$-simplicial set $\Sigma(A,\fg)$ as follows:

For $a:[p]\to[q], n\in{\Bbb N}$ one has
$$\Sigma(A,\fg)(a)_n=\MC(A_{np}\otimes\fg^q).$$

Define also a simplicial set $\sigma(A,\fg)$ to be the inverse limit
of $\Sigma(A,\fg)$ as a functor from $\CM$ to $\simpl$. 
We provide now a sufficient condition for a map $f:A\to B$ of bisimplicial
commutative dg algebras to induce an acyclic fibration $\Sigma(f,\fg)$
for any cosimplicial dg Lie algebra. According to~\Lem{afM->ss} this
implies that $\sigma(f,\fg)$ is also an acyclic fibration.

\subsubsection{} Any bisimplicial abelian group $A$ gives rise to a functor
$$A:\simpl\times\simpl\lra\Ab$$
which is uniquely described by the following properties:

--- $A(\Delta^m,\Delta^n)=A_{mn}$

--- $A$ carries direct limits over each one of the arguments to inverse limit.

We will identify bisimplicial abelian groups with the functors they define.

\subsubsection{}
\begin{defn}{bi-match} The matching space $M_{mn}(A)$ of a bisimplicial
abelian group $A$ is defined to be
$$ M_{mn}(A)=A(\partial\Delta^m,\Delta^n)\times
_{A(\partial\Delta^m,\partial\Delta^n)} A(\Delta^m,\partial\Delta^n)$$
where $\partial\Delta^n$ is the boundary of the $n$-simplex. 
One has a canonical map $A_{mn}\to M_{mn}(A).$
\end{defn}

\subsubsection{}
\begin{defn}{bi-af} A map $f:A\to B$ in $(\Delta^0)^2C(\Bbb Z)$ is called an
{\em acyclic fibration} if for any $m,n$ the canonical map
$$ A_{mn}\to B_{mn}\times_{M_{mn}(B)} M_{mn}(A)$$
is a surjective quasi-isomorphism.
\end{defn}

\subsubsection{}
\begin{prop}{bi-af->af}Let $f:A\to B$ in $(\Delta^0)^2\dgc(k)$ be an acyclic
fibration. Let $\fg$ be a cosimplicial nilpotent dg Lie algebra. Then
the induced map $\Sigma(f,\fg)$ is an acyclic fibration of $\CM$-simplicial
sets.
\end{prop}
\begin{pf}
This is a  direct calculation using~\ref{AFsur},~\ref{A+L}(a). Here we use 
that the natural map from $\fg^{n+1}$ to the $n$-th matching space $M^n(\fg)$ 
(see~\cite{bk}, X.5) is surjective.
\end{pf}
  
Now we wish to  formulate a sufficient condition for $f:A\to B$ to be an
acyclic fibration of bisimplicial dg algebras. The following trivial 
lemma will be useful.
\subsubsection{}
\begin{lem}{trivial}
Let in a commutative square  
$$
\begin{array}{ccc}
A & \overset{f}{\lra} & B \\
{{\scriptstyle g}\downarrow} &  & {{\scriptstyle h}\downarrow} \\
C & {\lra} & D \\
\end{array}
$$
of abelian groups the map $g$ and the map $\Ker(g)\to\Ker(h)$ be 
surjective. Then the induced map $A\to B\times_DC$ is also surjective.
\end{lem}

We start with a simplicial
case. Recall that a simplicial abelian group $A\in\Delta^0\Ab$ defines a
functor $A:\simpl\to\Ab$ by the formula
$$ A(S)=\Hom(S,A).$$

\subsubsection{}
\begin{lem}{simpl-af}Let $f:A\to B$ be a map in 
$\Delta^0C^{\geq 0}(\Bbb Z)=C^{\geq 0}(\Delta^0\Ab)$.
Suppose that

(a) $f_n:A_n\to B_n$ is a quasi-isomorphism in $C(\Bbb Z)$

(b) for any $S\in\simpl$ the map $f(S):A(S)\to B(S)$ is surjective

(c) for any $d\in\Bbb Z$ the $d$-components $A^d$ and $B^d$ are contractible 
simplicial abelian groups.

Then for any  injective map $\alpha:S\to T$ the induced map
$$ A(T)\to A(S)\times_{B(S)}B(T)$$
is a surjective quasi-isomorphism.
\end{lem}
\begin{pf}
For any $S\in\simpl$ the map $f(S):A(S)\to B(S)$ is quasi-isomorphism ---
this follows from~\cite{le}, remark at the end of III.2, applied to~\cite{ha},
Thm.~1.5.1.
This immediately implies that for any $\alpha:S\to T$ the induced map
$$ A(T)\to A(S)\times_{B(S)}B(T)$$
is a quasi-isomorphism. Let us prove that it is surjective if $\alpha$ is
injective.

Since $A^d$ are contractible (and Kan), the map $A(\alpha):A(T)\to A(S)$ is 
surjective. Thus, by~\Lem{trivial}, it suffices to check that the map
$$ \Ker A(\alpha)\to\Ker B(\alpha)$$ is surjective.

Since the functors $A,B:\simpl\to C(\Bbb Z)$ carry colimits to limits,
the map $\Ker A(\alpha)\to\Ker B(\alpha)$ is a direct summand of the map
$f(T/S): A(T/S)\to B(T/S)$ which is surjective by (b).
\end{pf}

Now we are able to prove the following criterion for a map of bisimplicial
complexes of abelian groups to be an acyclic fibration.

\subsubsection{}
\begin{prop}{criterion}
Let $f:A\to B$ be a map in $(\Delta^0)^2C^{\geq 0}(\Bbb Z)$ satisfying:

(a) for any $S,T\in\simpl$ the map $f(S,T)$ is surjective

(b) for any $S\in\simpl,\ p\in\Bbb N$ the map $f(S,\Delta^p):
A(S,\Delta^p)\to B(S,\Delta^p)$ is a quasi-isomorphism in $C(\Bbb Z)$.

(c) for any $S\in\simpl$ and any $d$ the simplicial abelian groups 
$A^d(S,\_)$, $B^d(S,\_)$ are contractible.

(d) for any $p,d$ the simplicial abelian group $A^d(\_,\Delta^p)$ is
 contractible.

Then $f$ is an acyclic fibration.
\end{prop}
\begin{pf}Apply~\Lem{simpl-af} to the map $f(S,\_):A(S,\_)\to B(S\_)$.
We immediately get that for any $S\in\simpl$ and any injective map 
$\alpha:T\to T'$ the induced map
$$ A(S,T')\to A(S,T)\times_{B(S,T)}B(S,T')$$
is a surjective quasi-isomorphism. Put $T'=\Delta^p,\ T=\partial\Delta^p.$
The map $A_{np}\to B_{np}\times_{M_{np}(B)}M_{np}(A)$ is then automatically
quasi-isomorphism, an we have only to check it is surjective.
According to (d) the map 
$ A(S',T')\to A(S,T')$
is surjective for $T'=\Delta^p$. 
Define
$$X(T)=\Ker(A(S',T)\to A(S,T)),\ Y(T)=\Ker(B(S',T)\to B(S,T)).$$

The groups $X(T)$ and $Y(T)$ are direct summands of $A(S'/S,T)$ and of 
$B(S'/S,T)$ respectively. Hence the map $X\to Y$ satisfies the hypotheses 
of~\Lem{simpl-af}. Therefore the map 
$$X(T')\to X(T)\times_{Y(T)}Y(T')$$
is surjective by~\Lem{trivial}.  Proposition is proven.
\end{pf}

\section{The main theorem}
\label{thm}

Now we came to the main result of the paper.

Let $\fg\in\Delta\dgl(k)$ be a cosimplicial nilpotent dg Lie 
$k$-algebra. Suppose that $\fg$ is {\em finitely dimensional in the 
cosimplicial sense}, i.e. that the normalization 
$$N^n(\fg)=\{x\in\fg^n|\sigma^i(x)=0\text{ for all } i\}$$
vanishes for sufficiently big $n$.

\subsection{}
\begin{thm}{main}
There is a natural homotopy equivalence
$$\Sigma(\Tot(\fg))\lra \Tot(\Sigma({\fg}))$$
in $\Kan$.
\end{thm}

Taking into account \ref{kan}, we easily get
\subsubsection{}
\begin{cor}{cor(main)}
Let $\fg\in\Delta\dgl(k)$ be a nilpotent cosimplicial dg Lie $k$-algebra.
 Suppose that $\fg$ is finitely dimensional in the 
cosimplicial sense.

Then there is a natural equivalence of groupoids
$$ \CC(\Tot(\fg))\ra\Tot(\CC(\fg)).$$
\end{cor}
\begin{pf}
One has to check that the functor $\Tot$ carries the map 
$\tau:\Sigma(\fg)\to\CC(\fg)$ to a homotopy equivalence. By~\cite{bk},~X.5,
it suffices to check $\tau$ is a fibration in sense of {\em loc. cit.}

This follows from (the second claim of)~\Prop{kan} since for any $n$
the natural map from $\fg^{n+1}$ to the $n$-th matching space 
$M^n(\fg)$ (in notations of {\em loc. cit}) is surjective.
\end{pf}

\subsection{Proof of the theorem}
The set of $n$-simplices of the left-hand side is
$$ \MC(\Omega_n\otimes\invlim \Omega_p\otimes\fg^q)$$
and for the right-hand side:
$$ \invlim\MC(\Omega(\Delta^n\times\Delta^p)\otimes\fg^q).$$
Here the inverse limits are taken over the category $\CM$ defined
in~\ref{catM} and~\Lem{repr} is used to get the second
formula.

Taking into account that the functor $\MC$ commutes with the inverse limits,
a canonical map $\Sigma\circ{\Tot(\fg)}\to\Tot\circ\Sigma({\fg})$ is 
defined by the composition
$$\Omega_n\otimes\invlim \Omega_p\otimes\fg^q\to
\invlim \Omega_n\otimes\Omega_p\otimes\fg^q\to
\invlim\Omega(\Delta^n\times\Delta^p)\otimes\fg^q$$
the latter arrow being induced by the canonical projections of
$\Delta^n\times\Delta^p$ to $\Delta^n$ and to $\Delta^p$.

We wish to prove that the map described induces  a homotopy equivalence.

Since $\fg$ is finitely dimensional, the first map in the composition is 
bijective ---   see~\cite{hdtc}, Thm. 6.11. In order to prove that the 
second map is a homotopy equivalence, let us fix $n$ and $p$ and 
present the map 
$$\alpha:\Omega_n\otimes\Omega_p\to\Omega(\Delta^n\times\Delta^p)$$
as the composition
$$\Omega_n\otimes\Omega_p\overset{\beta}{\lra}
\Omega(\Delta^n\times\Delta^{p+1})\otimes\Omega_p\overset{\pi}{\lra}
\Omega(\Delta^n\times\Delta^p)$$
where $\beta$ is induced by the projection 
$\Delta^n\times\Delta^{p+1}\to\Delta^n$
and $\pi$ by the pair of maps $\id\times\partial^{p+1}:\Delta^n\times
\Delta^p\to\Delta^n\times\Delta^{p+1}$, 
$\pr_2:\Delta^n\times\Delta^p\to\Delta^p$.

We will check below that the maps $\beta$ and $\pi$ induce homotopy
equivalences for different reasons: $\beta$ induces a strong deformation
retract and $\pi$ induces an acyclic (Kan) fibration. This will prove
the theorem.

\subsubsection{Notations} 

Define $\Delta^{+1}$ to be the cosimplicial simplicial set with
$(\Delta^{+1})^n=\Delta^{n+1}$ whose cofaces and codegeneracies are the
standard maps between the standard simplices (they all preserve the final
vertex).

Put $A_{np}=\Omega_n\otimes\Omega_p$, $B_{np}=\Omega(\Delta^n\times
\Delta^p)$, $C_{np}=\Omega(\Delta^n\times\Delta^{p+1})\otimes\Omega_p$.
$A$ and $B$ are bisimplicial commutative dg algebras by an obvious reason.
Bisimplicial structure on $C$ in defined by the cosimplicial structure
on $\Delta^{+1}$.

Our aim is to prove that the maps $\beta:A\to C$ and $\pi:C\to B$
induce homotopy equivalences $\sigma(\beta,\fg)$ and $\sigma(\pi,\fg)$.

We will check immediately that $\sigma(\beta,\fg)$ is a strong 
deformation retract. Afterwards, using the criterion~\ref{criterion}
we will get that $\sigma(\pi,\fg)$ is an acyclic fibration.

\subsubsection{Checking $\beta$} 

For $S\in\simpl$ and nilpotent $\fg\in\dgl(k)$ denote by 
$\Sigma^S(\fg)$ (or just $\Sigma^S$ when $\fg$ is one and the same)
the simplicial set whose set of $n$-simplicies is 
$\MC(\Omega(\Delta^n\times S)\otimes\fg)$.

Any map $f:K\times S\to T$ in $\simpl$ induces a map
$$\Sigma^f:K\times \Sigma^T\to \Sigma^S$$ 
as follows. Let 
$k\in K_n,\ x\in\Sigma^T_n=\MC(\Omega(\Delta^n\times T)\otimes\fg)$. 
Denote by $F_k$ the composition
$$ \Delta^n\times S\overset{\diag\times 1}{\lra}
\Delta^n\times\Delta^n\times S\overset{1\times k\times 1}{\lra}
\Delta^n\times K\times S\overset{1\times f}{\lra}
\Delta^n\times T.$$
Then $\Sigma^f(k,x)$ is defined to be $(\Omega(F_k)\otimes\id_{\fg})(x)$.

Define the map $\Phi_p:\Delta^1\times
\Delta^{p+1}\to\Delta^{p+1}$ as the one given on the level of posets
by the formula
$$ \Phi_p(i,j)=\begin{cases}
                j \text{ if } i=0\\
                p+1\text{ if } i=1
              \end{cases}
$$

The maps $\Phi_p$ induce the maps 
$\Sigma^{\Phi_p}(\Omega_p\otimes\fg^q)$ which define for any
$a:[p]\to [q]$ in $\CM$ the simplicial set 
$\Sigma(A,\fg)(a)=\MC(A_{\bullet p}\otimes\fg^q)$ as a strong
deformation retract of $\Sigma(C,\fg)(a)$. 
The retractions $\Sigma^{\Phi_p}(\Omega_p\otimes\fg^q)$ are functorial in 
$a\in\CM$, therefore the inverse image map $\sigma(\beta,\fg)$ is also strong
deformation retract.

\subsubsection{Checking $\pi$} Now we will prove that the map
$\pi:C\to B$ of bisimplicial dg algebras defined above induces an acyclic
fibration $\sigma(\pi,\fg)$. 

Let us check the hypotheses of~\ref{criterion}. 

First of all, let us check that $\pi(S,T)$ is surjective. For this
consider $D_{np}=\Omega(\Delta^n,\Delta^{p+1})$ and the map
$\rho_{np}:D_{np}\to B_{np}$ which is the composition of $\pi$
with the natural embedding $D_{np}\to C_{np}=D_{np}\otimes\Omega_p$.
Surely, it suffices to prove that $\rho(S,T)$ is surjective.
For this we note that $D(S,T)=\Omega(S\times\cone{T})$ where
$\cone{\ }:\simpl\to\simpl$ is the functor satisfying the condition
$\cone{\ }|_{\Delta}=\Delta^{+1}$ and preserving colimits 
(see~\cite{gz},II.1.3).  Then, since the map 
$S\times T\to S\times\cone{T}$ is injective, the map
$$\Omega(S\times\cone{T})\to \Omega(S\times T)$$
is surjective.

Next, one has 
$$C(S,\Delta^p)=\Omega(S\times\Delta^{p+1})\otimes\Omega_p,$$
$$B(S,\Delta^p)=\Omega (S\times\Delta^p)$$
so the condition~\ref{criterion} (b) is fulfilled.

The simplicial abelian groups $C^d(S,\_)$ are contractible by the K\"unneth
formula.  The abelian groups $B^d(S,\_)=\Omega^d(S\times\_)$
are contractible since any injective map $T\to T'$ induces a
surjection $\Omega(S\times T')\to\Omega(S\times T)$. The same reason
proves the condition (d).

Therefore, the map $\pi:C\to B$ of bisimplicial commutative dg algebras 
is an acyclic fibration by~\Prop{criterion} and then by~\Prop{bi-af->af}
the map $\sigma(\pi,\fg)$ is an acyclic fibration in $\simpl$.

The Theorem is proven.

\section{Application to Deformation theory}
\label{app}

In this Section we describe how to deduce from~\Cor{cor(main)}
the description of universal formal deformations for some typical
deformation problems.

Consider, for example, three deformation problems which have been
studied in~\cite{hdtc}.

 Let $X$ be a smooth separated scheme $X$ over a field $k$ 
of characteristic $0$, $G$ an algebraic group over $k$ and $p: P\lra X$ 
a $G$-torsor over $X$. Consider the following deformation 
problems.    

{\bf Problem 1.} Flat deformations of $X$. 

{\bf Problem 2.} Flat deformations of the pair $(X,P)$. 

{\bf Problem 3.} Deformations of $P$ ($X$ being fixed).

To each problem one can assign a sheaf of $k$-Lie algebras $\fg_i$ on $X$
(these are the sheaves ${\cal A}_i,\ i=1,2,3$ from {\em loc. cit.},
Section 8).

According to Grothendieck, to each problem corresponds a (2-)functor 
of infinitesimal deformations
$$
F_i:\art/k\lra\Grp
$$
from the category of local artinian $k$-algebras with the residue 
field $k$ to the (2-)category of groupoids.  

In each case, $\fg_i$ is "a sheaf of infinitesimal automorphisms" 
corresponding to $F_i$ (in the sense of ~\cite{sga1}, Exp.III, 5,
especially Cor. 5.2 for Problem 1; for the other problems the meaning is 
analogous).

If $X$ is affine, the functor $F_i$ is equivalent to the Deligne groupoid
$\CC_L,\ L=\Gamma(X,\fg_i)$ defined as the functor
$$ \CC_L:\art/k\to\Grp$$
which is given by the formula
$$\CC_L(A)=\CC(\fm\otimes L)$$
where $\fm$ is the maximal ideal of $A$.

The deformation functor $F_i$ defines a stack $\CF_i$ in the Zariski topology
 of $X$. The Deligne functor defines a fibered category $\CC_{\fg_i}$ which 
assigns a groupoid $\CC_{\Gamma(U,\fg_i)}(A)$ to each Zariski open set $U$ and 
to each $A\in\art/k$.

We have a canonical map of fibered categories $ \CC_{\fg_i}\to\CF_i$ so that
$\CF_i$ is equivalent to the stack associated to the fibered category 
$\CC_{\fg_i}$.

Using~\Cor{cor(main)} we immediately get

\subsection{}
\begin{cor}{}
The deformation groupoid $F_i$ is naturally equivalent to the Deligne 
groupoid associated with the dg Lie algebra $\Right\Gamma^{\Lie}(X,\fg_i)$.
\end{cor}

In particular, the following generalization of~\cite{hdtc}, Thm.~8.3
takes place.

\subsection{}
\begin{cor}{generalization}
Suppose that $H^0(X,\fg_i)=0$; let $\fS=\Spf(R)$ be the base 
of the universal formal deformation for Problem $i$.
Then we have a canonical isomorphism
$$
R^*=H^{Lie}_0(\Right\Gamma^{Lie}(X,\fg_i))
$$
where $R^*$ denotes the space of continuous $k$-linear maps $R\lra k$ 
($k$ considered in the discrete topology) and $H^{\Lie}_0$ denotes the 
$0$-th Lie homology.
\end{cor}

Recall that this result has been proven  in~{\em loc. cit.} for $\fS$
formally smooth.

\end{document}